%% file: LiuEtal11ITA.tex
\documentclass[10pt,conference]{IEEEtran}
\usepackage{times,amsmath,amssymb,graphicx,epsf,psfrag,epsfig,cite}

\input macros

\def\ie{{\it i.e.,\ \/}}

\def\defeq{{\,\stackrel{\Delta}{=}}\,}

\def\scalefig#1{\epsfxsize #1\textwidth}

\newcommand{\mbbE}{\mathbb{E}}

\newtheorem{theorem}{Theorem}

\begin{document}

\title{Decentralized Restless Bandit with Multiple Players and Unknown Dynamics}

\author{\authorblockN{Haoyang Liu, ~~~  Keqin Liu, ~~~ Qing Zhao}
\authorblockA{Department of Electrical and Computer Engineering\\
University of California, Davis, CA 95616\\
 \{liu, kqliu, qzhao\}@ucdavis.edu}}

\maketitle

\begin{abstract}
We\footnotetext{0This work was supported by the Army Research Office under Grant W911NF-08-1-0467 and by the National Science
Foundation under Grant CCF-0830685.}
consider decentralized restless multi-armed bandit problems with
unknown dynamics and multiple players. The reward state of each arm
transits according to an unknown Markovian rule when it is played
and evolves according to an arbitrary unknown random process when it
is passive. Players activating the same arm at the same time collide
and suffer from reward loss. The objective is to maximize the
long-term reward by designing a decentralized arm selection policy
to address unknown reward models and collisions among players. A
decentralized policy is constructed that achieves a regret with
logarithmic order when an arbitrary nontrivial bound on certain system parameters is known.
When no knowledge about the system is available, we extend the policy to achieve a regret
arbitrarily close to the logarithmic order. The result finds applications in communication
networks, financial investment, and industrial engineering.

\end{abstract}

\section{Introduction}\label{sec:introduction}

\subsection{The Classic MAB with A Single Player}
In the classic MAB, there are $N$ independent arms and a single
player. Each arm, when played, offers an i.i.d. random reward to the
player. The reward distribution of each arm is unknown. At each
time, the player chooses one arm to play, aiming to maximize the
total expected reward in the long run. This problem involves the
well-known tradeoff between exploitation and exploration. For
exploitation, the player should select the arm with the largest
sample mean of reward. For exploration, the player should select an
under-played arm to learn its reward statistics.

Under the non-Bayesian formulation, the performance measure of an
arm selection policy is the so-called \emph{regret} or \emph{the
cost of learning} defined as the reward loss with respect to the
case with known reward models~\cite{Lai&Robbins85AAM}. In~1985 Lai
and Robbins showed that the minimum regret grows at a logarithmic
order under certain regularity conditions~\cite{Lai&Robbins85AAM}.
The best leading constant was also obtained, and an optimal policy
was constructed to achieve the minimum regret growth rate (both the
logarithmic order and the best leading constant). In 1987,
Anantharam~\etal extended Lai and Robbins's results to accommodate
multiple simultaneous plays~\cite{Anantharam:87-1} and a Markovian
reward model where the reward of each arm evolves as an unknown
Markov process over successive plays and remains frozen when the arm
is passive (the so-called \emph{rested} Markovian reward
model)~\cite{Anantharam:87-2}.

Several other simpler policies have been developed to achieve
logarithmic regret for the classic MAB under an i.i.d. reward
model~\cite{Agrawal:95,Auer:02}. In particular, the index
policy---referred to as Upper Confidence Bound 1 (UCB-1)---proposed
in~\cite{Auer:02} achieves the logarithmic regret with a uniform
bound on the leading constant over time. In~\cite{Tekin:10}, UCB-1
was extended to the rested Markovian reward model adopted
in~\cite{Anantharam:87-2}.

\subsection{Decentralized MAB with Distributed Multiple Players}
In~\cite{Keqin:10}, Liu and Zhao formulated and studied a
decentralized version of the classic MAB with $M$ ($M<N$)
distributed players under the i.i.d. reward model. Different arms
can have different reward distributions and they are unknown to the
players. At each time, a player chooses one arm to play based on its
{\em local} observation and decision history without exchanging
information with other players. Collisions occur when multiple
players choose the same arm, and, depending on the collision model,
either no one receives reward or the colliding players share the
reward in an arbitrary way. The objective is to maximize the
long-term sum reward from all players. Another desired feature of
policies for decentralized MAB is fairness, \ie different players
have the same expected reward growth rate. Liu and Zhao proposed the
Time Division Fair Sharing (TDFS) framework,  it achieves the same
logarithmic regret order as the centralized case where all players
share their observations in learning and collisions are eliminated
through centralized perfect scheduling~\cite{Keqin:10}. Assuming a
Bernoulli reward model, decentralized MAB was also addressed in
\cite{Anima}, where the single-player policy UCB-1 was extended to
the multi-player setting.

\subsection{Main Results}
In this paper, we consider the decentralized MAB with a restless
Markovian reward model. In a single-player restless MAB, the reward
state of each arm transits according to an unknown Markovian rule
when played and transits according to an arbitrary unknown random
process when passive as addressed in our prior
work~\cite{Liu&Liu&Zhao10ICASSP}. In~\cite{Liu&Liu&Zhao10ICASSP}, we
proposed a policy Restless UCB (RUCB), which achieves a logarithmic
order of the weak regret defined as the reward loss compared to the
case when the player knows which arm is the most rewarding and
always plays the best arm. RUCB borrows the index form of UCB-1
given in~\cite{Auer:02} and has a deterministic epoch structure with
carefully chosen epoch lengths to balance exploration and
exploitation. The concept of weak regret was first used in
\cite{Auer-nonsto}; it measures the reward loss with respect to the
optimal \emph{single-arm} policy, which, while optimal under the
i.i.d. and rested Markovian reward models (up to an $O(1)$ term of loss for the latter),
is no longer optimal in
general under a known restless reward model. Analysis of the strict
regret of restless MAB is in general intractable given that finding
the optimal policy of a restless bandit under \emph{known} model is
itself PSPACE-hard in general~\cite{Papadimitriou:99}.

In this paper, we extend RUCB proposed in our prior
work~\cite{Liu&Liu&Zhao10ICASSP} to a decentralized setting of
restless MAB with multiple players. We consider two types of
restless reward models: exogenous restless model and endogenous
restless model. In the former, the system itself is rested: the
state of an arm does not change when the arm is not engaged.
However, from each individual player's perspective, arms are
restless due to actions of other players that are unobservable and
uncontrollable. Under the endogenous restless model, the state of an
arm evolves according to an arbitrary unknown random process even
when the arm is not played. Under both restless models, we extend
RUCB to achieve a logarithmic order of the regret. The result for
the exogenous restless model, however, is stronger in the sense that
the regret is indeed defined with respect to the optimal policy
under known reward models. This is possible due to the inherent
\emph{rested} nature of the systems.

There are a couple of parallel work to~\cite{Liu&Liu&Zhao10ICASSP}
on the single-player restless MAB. In~\cite{Tekin:10-2}, Tekin and
Liu adopted the weak regret and proposed a policy that achieves
logarithmic (weak) regret when certain knowledge about the system
parameters is available~\cite{Tekin:10-2}. The policy proposed
in~\cite{Tekin:10-2} also uses the index form of UCB-1 given
in~\cite{Auer:02}, but the structure is different from RUCB proposed
in~\cite{Liu&Liu&Zhao10ICASSP}. Specifically, under the policy proposed
in~\cite{Tekin:10-2}, an arm is played consecutively for a random number
of times determined by the regenerative cycle of a particular state, and
observations obtained outside the regenerative cycle are not used in learning.
RUCB, however, has a deterministic epoch structure, and all observations
are used in learning. In~\cite{QZ10}, the strict regret
was considered for a special class of restless MAB. Specifically,
when arms are governed by stochastically identical two-state Markov
chains, a policy was constructed in~\cite{QZ10} to achieve a regret
with an order arbitrarily close to logarithmic.

\noindent{\bf Notation} For two positive integers $k$ and $l$,
define $k\oslash l\defeq ((k-1)~\mbox{mod}~l)+1$, which is an
integer taking values from $1,2,\cdots,l$.

\section{Problem Formulation}\label{sec:problemformulation}
In the decentralized MAB problem, we have $M$ players and $N$
independent arms. At each time, each player chooses one arm to play.
Each arm, when played (activated), offers certain amount of reward
that models the current state of the arm. Let $s_j(t)$  and
$\mathcal{S}_j$ denote the state of arm $j$ at time $t$ and the
state space of arm $j$ respectively. Different arms can have
different state spaces. When arm $j$ is played, its state changes
according to a Markovian rule with $P_j$ as the transition matrix.
The transition matrixes are assumed to be irreducible, aperiodic,
and reversible. As for the state transition of passive arms, we
consider two models: endogenous restless model and exogenous
restless model. In the endogenous restless model, arm states change
in arbitrary ways when not played. In the exogenous restless model,
arm states remain frozen if not engaged. The players do not know the
transition matrices of the arms and do not communicate with each
other. Conflicts occur when different players select the same arm to
play. Under different conflict models, either the players in
conflict share the reward or no one obtains any reward. The
objective is to maximize the expected total reward collected in the
long run. Let $\vec{\pi_j}=\{\pi_s^j\}_{s\in\mathcal{S}_j}$ denote
the stationary distribution of arm $j$ (under $P_j$), where
$\pi_s^i$ is the stationary probability (under $P_j$) that arm $j$
is in state $s$. The stationary mean reward $\mu_j$ is given by $
\mu_j= \sum_{s \in \mathcal{S}_j} s \pi_s^j$. Let $\sigma$ be a
permutation of $\{1, \cdots, N\}$ such that
\[
\mu_{\sigma(1)} \geq \mu_{\sigma(2)} \geq \mu_{\sigma(3)} \geq
\cdots \geq \mu_{\sigma(N)}.
\]
A policy $\Phi$ is a rule that specifies an arm to play based on the
local observation history. Let $t_j(n)$ denote the time index of the
$n$th play on arm $j$, and $T_j(t)$ the total number of plays on arm
$j$ by time $t$. Notice that both $t_j(n)$ and $T_j(t)$ are random
variables with distributions determined by the policy $\Phi$. Under
the conflict model where players in conflict share the reward, the
total reward by time $t$ is given by
\begin{eqnarray}
R(t) = \sum_{j=1}^{N}\sum_{n=1}^{T_j(t)}s_j(t_j(n)).
\end{eqnarray}

Under the conflict model where no players in conflict obtain any
reward, the total reward by time $t$ is given by
\begin{eqnarray}
R(t) = \sum_{j=1}^{N}\sum_{n=1}^{T_j(t)}s_j(t_j(n))
\mathbb{I}_j({t_j(n)}).
\end{eqnarray}
where $\mathbb{I}_j({t_j(n)}) = 1$ if arm $j$ is played by one and
only one player at time $t_j(n)$, and $\mathbb{I}_j({t_j(n)}) = 0$
otherwise.

As mentioned in Sec.~\ref{sec:introduction}, for both restless
models, performance of any policy $\Phi$ is evaluated using regret
$r_{\Phi} (t)$ defined as the reward loss with respect to having $M$
best arms constantly engaged. Specifically, for both restless
models, regret is defined as follows:
\begin{eqnarray} \label{eqn:regret}
r_\Phi(t) = t\sum_{i=1}^{M}\mu_{\sigma(i)}- \mbbE_{\Phi}R(t)+ O(1),
\end{eqnarray}
where the constant $O(1)$ is caused by the transient effects of
playing the $M$ best arms, $\mbbE_{\Phi}$
denotes the expectation with respect to the random process induced
by policy $\Phi$. The objective is to minimize the growth rate of
the regret. Note that the constant $O(1)$ term can be ignored when
studying the growth rate of the regret.

\begin{figure}[h]
\centerline{
\begin{psfrags}
\psfrag{t=}[c]{\small Slot} \psfrag{Exploitation epochs}[l]{\small
Exploitation epochs} \psfrag{Exploration epochs}[l]{\small
Exploration epochs}\psfrag{t=}[c]{\small Slot}
\psfrag{Epoch}[c]{\small Epoch} \psfrag{arm}[l]{\scriptsize arm}
\psfrag{Policy structure}[c]{\small The general structure of
decentralized RUCB}
\psfrag{1}[c]{\scriptsize$1$}\psfrag{2}[c]{\scriptsize
$2$}\psfrag{3}[c]{\scriptsize $3$}\psfrag{4}[c]{\scriptsize $4$}
\psfrag{5}[c]{\scriptsize$5$}\psfrag{6}[c]{\scriptsize$6$}
\psfrag{7}[c]{\scriptsize$7$} \psfrag{8}[c]{\scriptsize$8$}
\psfrag{2^(n-1)}[c]{\scriptsize$2^{n-1}$}
\psfrag{2*2^(n-1)}[r]{\scriptsize$ 2 \times 4^{n-1}$}
\psfrag{12}[l]{\scriptsize$2 \times 4^{n-1}$}
\psfrag{123}[l]{\tiny$2(M-1) \times 4^{n-1}$}
\psfrag{1234}[l]{\scriptsize$2M \times 4^{n-1}$}
\psfrag{a*1}[l]{\scriptsize$a^*_2$} \psfrag{a*N}[l]{\tiny
$a^*_{M-1}$} \psfrag{a*M}[l]{\scriptsize$a^*_{M}$}
\psfrag{(N-1)*2^(n-1)+1}[c]{\tiny$(N-1) \times 4^{n-1}+1$}
\psfrag{N*2^(n-1)}[l]{\scriptsize$N \times 4^{n-1}$}
\psfrag{N}[c]{\scriptsize$N$}
 \psfrag{Structure of exploration epoch}[c]{\small Structure of the $n$th exploration epoch for player $1$}
  \psfrag{Structure of exploitation epoch}[c]{\small Structure of the $n$th exploitation epoch  for player $1$}
\psfrag{Comment for first time}[l]{\scriptsize Compute the indexes
and identify the arms with the $M$ highest indexes}
\psfrag{a*}[l]{\scriptsize $a^*_1$} \psfrag{First}[l]{\scriptsize
First} \psfrag{Second}[l]{\scriptsize
Second}\psfrag{Third}[l]{\scriptsize Third} \psfrag{Two
types}[l]{\scriptsize Two types of epochs are not always
interleaving.}\psfrag{exploration}[l]{\scriptsize exploration}
\psfrag{exploitation}[cl]{\scriptsize exploitation}
\psfrag{epoch}[l]{\scriptsize epoch} \psfrag{ Policy
structure}[l]{\small Policy structure}\psfrag{...}[c]{\Large
$\cdots$} \scalefig{0.48}\epsfbox{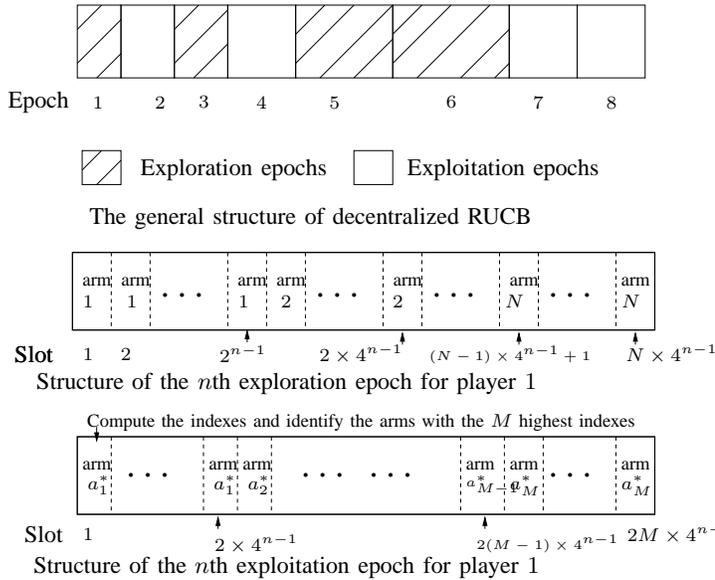} \psfrag{Two
types}[c]{\scriptsize Two types of epochs are not always
interleaving.}
\end{psfrags}}
\caption{Epoch structures of decentralized RUCB}
\label{fig:structure}
\end{figure}

\section{The Decentralized RUCB Policy}\label{sec:policy}

The proposed decentralized RUCB is based on an epoch structure. We
divide the time into disjoint epochs. There are two types of epochs:
exploitation epochs and exploration epochs (see an illustration in
Fig.~\ref{fig:structure}). In the exploitation epochs, the players
calculate the indexes of all arms and play the arms with the $M$
highest indexes, which are believed to be the $M$ best arms. In the
exploration epochs, the players obtain information of all arms by
playing them equally many times. The purpose of the exploration
epochs is to make decisions in the exploitation epochs sufficiently
accurate. As shown in Fig. \ref{fig:structure}, in the $n$th
exploration epoch, each player plays every arm $4^{n-1}$ times. At
the beginning of the $n$th exploitation epoch the player calculates
index for every arm (see~(\ref{eqn:index}) in
Fig.~\ref{fig:decenpolicy}) and selects the arm with the $M$ highest
indexes (denoted as arm $a^*_{(1)}$ to arm $a^*_{(M)}$). Each
exploitation epoch is divided into $M$ subepochs with each having a
length of $2 \times 4^{n-1}$. Player $k$ plays arm $a^*_{((m-k+M+1)
\oslash M)}$ in the $m$th subepoch of each exploitation epoch. The
details on interleaving the two types of epochs are given in Step
$2$ in Fig.~\ref{fig:decenpolicy}. Specifically, whenever
sufficiently many ($D\ln t$, see~(\ref{eqn:condition})) observations
have been obtained from every arm in the exploration epochs, the
player is ready to proceed with a new exploitation epoch. Otherwise,
another exploration epoch is required to gain more information about
each arm. It is also implied in (\ref{eqn:condition}) that only
logarithmically many plays are spent in the exploration epochs,
which is one of the key reasons for the logarithmic regret of
decentralized RUCB. This also implies that the exploration epochs
are much less frequent than the exploitation epochs. Though the
exploration epochs can be understood as the ``information
gathering'' phase, and the exploitation epochs as the ``information
utilization'' phase, observations obtained in the exploitation
epochs are also used in learning the arm dynamics. This can be seen
in Step $3$ in Fig.~\ref{fig:decenpolicy}. The epoch structure, (\ie
the starting and ending points of epochs) are prefixed numbers only
depending on parameter $D$. This is one of the key reasons why
different players can be coordinated (\ie entering the same epoch at
the same time) without intercommunications.

\begin{figure}[htb]
\begin{center}
\noindent\fbox{
\parbox{3.4in}
{ \centerline{\underline{{\bf Decentralized RUCB}}} {\small Time is
divided into epochs. There are two types of epochs, exploration
epochs and exploitation epochs. At the beginning of the $n$th
exploitation epoch, we choose the $M$ arms to play, each of them for
$2 \times 4^{n-1}$ many times. In the $n$th exploration epoch, we
play every arm $4^{n-1}$ many times. Let $n_O(t)$ denote the number
of exploration epochs played by time $t$ and $n_I(t)$ the number of
exploitation epochs played by time $t$.
\begin{enumerate}
\item[1.] At $t=1$, we start the first exploration epoch, in which every arm is played once. We set $n_O(N+1) = 1$, $n_I(N+1) = 0$. Then go to Step $2$.
\item[2.] Let $X_1(t) = (4^{n_O(t)}-1)/3$ be the time spent on each arm in exploration epochs by time
$t$. Choose $D$ according to
(\ref{eqn:conditionL})(\ref{eqn:conditionD}). If
\begin{eqnarray}\label{eqn:condition}
X_1(t)
> D \ln t ,
\end{eqnarray} go to Step $3$ (start an exploitation epoch). Otherwise, go to
Step $4$ (start an exploration epoch).
\item[3.] Calculate indexes $d_{i,t}$ for all arms using the formula below:
\begin{eqnarray}\label{eqn:index}
d_{i,t} = \bar{s}_i(t) + \sqrt{\frac{L \ln t}{T_i(t)}},
\end{eqnarray} where $t$ is the current
time, $\bar{s}_i(t)$ is the sample mean from arm $i$ by time $t$,
$L$ is chosen according to (\ref{eqn:conditionL}), and $T_i(t)$ is
the number of times we have played arm $i$ by time $t$. Then choose
the arms with the $M$ highest indexes (arm $a^*_{(1)}$ to arm
$a^*_{(M)}$). Each exploitation epoch is divided into $M$ subepochs
with each having a length of $2 \times 4^{n-1}$. Player $k$ plays
arm $a^*_{((m-k+M+1 )\oslash M)}$ in the $m$th subepoch of each
exploitation epoch. After arm $a^*_{(1)}$ to arm $a^*_{(M)}$ are
played, increase $n_I$ by one and go to step $2$.
\item[4.] Each Play each arm for $4^{(n_O-1)}$ slots. Each exploration epoch is divided into $N$ subepochs
with each having a length of $4^{(n_O-1)}$. Player $k$ plays arm
$a^*_{(m-k+N+1 \oslash N)}$ in the $m$th subepoch of each
exploitation epoch. After all the arms are played, increase $n_I$ by
one and go to step $2$.
\end{enumerate}}
}} \caption{Decentralized RUCB policy}\label{fig:decenpolicy}
\end{center}
\end{figure}

\subsection{Eliminate Pre-Agreement}

So far we have assumed a pre-agreement among the players: they
target at the $M$ best arms with different offsets to avoid
excessive collisions. In this subsection, we show that this
pre-agreement can be eliminated while maintaining the logarithmic
order of the system regret. Furthermore, players can join the system
at different times without any global synchronization. Specifically,
at each player, the structure of the exploration and exploitation
epochs is the same as the local RUCB policy with pre-agreement. The
only difference here is that in each exploitation epoch, the player
randomly chooses one of the $M$ arms considered as the best to play
whenever a collision with other players is observed. If no collision
is observed, the player keeps playing the same arm. This simple
elimination of pre-agreement leads to a complete decentralization
among players while achieving the same logarithmic order of the
system regret. Except that each player can join the system according
to the local schedule, the player can also leave the system for an
arbitrary finite time period.

\section{The Logarithmic Regret of decentralized RUCB} \label{sec:Performance}
In this section, we show that the regret achieved by the
decentralized RUCB policy has a logarithmic order. This is given in
the following theorem.
\begin{theorem}\label{thm:singleuser}
Under the exogenous restless Markovian reward model, assume that
when arms are engaged, they can be modeled as finite state,
irreducible, aperiodic, and reversible Markov chains. All the states
(rewards) are positive. Let $\pi_{\min} = \min_{s\in \mathcal{S}_i,
1\leq i \leq N} \pi_s^i$, $\epsilon_{\max} = \max_{1\leq i \leq N}
\epsilon_i$, $\epsilon_{\min} = \min_{1\leq i \leq N} \epsilon_i$,
$s_{\max} = \max_{s\in \mathcal{S}_i, 1\leq i \leq N}s$, $s_{\min} =
\min_{s\in \mathcal{S}_i, 1\leq i \leq N}s$, and
$|\mathcal{S}|_{\max} = \max_{1\leq i \leq N} |\mathcal{S}_i|$ where
$\epsilon_i = 1- \lambda_{i}$ ($\lambda_{i}$ is the second largest
eigenvalue of the matrix $P_i$). Assume that different arms have
different $\mu$ values \footnote{This assumption can be relaxed by
utilizing the shared index set. This assumption is only for
simplicity of the presentation.} Set the policy parameters $L$ and
$D$ to satisfy the following conditions:
\begin{eqnarray}\label{eqn:conditionL}
L & \geq &\frac{1}{\epsilon_{\min}}(4 \frac {20s^2_{\max}
|\mathcal{S}|_{\max}^2}{(3-2\sqrt{2})}+ 10 s^2_{\max}),
\end{eqnarray}
\begin{eqnarray}\label{eqn:conditionD}
D & \geq & \frac{4L}{(\min_{j\leq M}( \mu_{\sigma(j)} -
\mu_{\sigma(j+1)}))^2}.
\end{eqnarray}

Under the conflict model where players share the reward, the regret
of decentralized RUCB at the end of any epoch can be upper bounded
by
\begin{eqnarray} \label{eqn:con1}
r_{\Phi}(t)& \leq &\frac{1}{3}[4(3D \ln t +1)-1]  \left(
\sum_{i=1}^{M}\mu_{\sigma(i)}
- \frac{M}{N} \sum_{i=1}^{N}\mu_{\sigma(i)}\right)\nonumber \\
 &  & + 3\lceil \log_4 (\frac{3}{2} (t- N) +1)
\rceil(1+ \frac{\epsilon_{\max} \sqrt{L}}{10s_{\min}}) \nonumber \\
 &  & \qquad \sum_{i=1}^{M-1} \sum_{j=1,j \neq i}^N
\mu_{\sigma(i)}\frac{|\mathcal{S}_{\sigma(i)}| + |\mathcal{S}_{\sigma(j)}|}{\pi_{\min}} \nonumber \\
&  & +  3 \lceil \log_4 (\frac{3}{2} (t- N) +1) \rceil(1+
\frac{\epsilon_{\max} \sqrt{L}}{10s_{\min}}) \nonumber \\
 &  & \qquad \sum_{j=M+1}^N
(\mu_{\sigma(M)}-\mu_{\sigma(j)})\frac{|\mathcal{S}_{\sigma(i)}| +
|\mathcal{S}_{\sigma(j)}|}{\pi_{\min}}. \nonumber \\&  & + 3 \lceil
\log_4 (\frac{3}{2} (t- N) +1) \rceil(1+ \frac{\epsilon_{\max}
\sqrt{L}}{10s_{\min}}) \nonumber \\
 &  & \qquad  \sum_{j=1}^{M-1} \mu_{\sigma(M)}
\frac{|\mathcal{S}_{\sigma(M)}| +
|\mathcal{S}_{\sigma(j)}|}{\pi_{\min}}\nonumber \\
&&+ \sum_{i=1}^{N}[(\min_{s\in \mathcal{S}_i} \pi_s)^{-1}
\sum_{s\in \mathcal{S}_i} s]
\end{eqnarray}

Under the model where no player in conflict gets any reward, the
regret of decentralized RUCB at the end of any epoch can be upper
bounded by:
\begin{eqnarray}\label{eqn:con2}
r_{\Phi}(t)& \leq &3 \lceil \log_4 (\frac{3}{2} (t- N) +1) \rceil(1+
\frac{\epsilon_{\max}
\sqrt{L}}{10s_{\min}})\nonumber \\
 &  & \qquad  (\sum_{i=1}^{M}\mu_{\sigma(i)}) \sum_{i=1}^{M}
\sum_{j=1,j \neq i}^N \frac{|\mathcal{S}_{\sigma(i)}| +
|\mathcal{S}_{\sigma(i)}|}{\pi_{\min}}\nonumber \\
 &  & +
\frac{1}{3}[4(3D \ln t +1)-1] \nonumber \\
 &  & \qquad \left( \sum_{i=1}^{M}\mu_{\sigma(i)}
- \frac{M}{N} \sum_{i=1}^{N}\mu_{\sigma(i)}\right)
\nonumber \\
&& + \sum_{i=1}^{N}[(\min_{s\in \mathcal{S}_i} \pi_s)^{-1}
\sum_{s\in \mathcal{S}_i} s]
\end{eqnarray}

\end{theorem}

We point out that upper bounds of regret in Theorem~$1$ can be
extended to any time $t$ instead of only for ending points of
epochs. They can also be extended to the endogenous restless model
in terms of weak regret. The no pre-agreement version of
decentralized RUCB can also achieve regret with a logarithmic order.
\begin{proof}
See Appendix~A for details.
\end{proof}

Theorem~$1$ requires an arbitrary (nontrivial) bound on
$s^2_{\max}$, $|\mathcal{S}|_{\max}$, $\epsilon_{\min}$, and
$\min_{j\leq M}( \mu_{\sigma(j)} - \mu_{\sigma(j+1)})$. In the case
where these bounds are unavailable, $D$ and $L$ can be chosen to
increase with time to achieve a regret order arbitrarily close to
logarithmic order. This is formally stated in the following theorem.
\begin{theorem}\label{thm:arbitraryclose}
Assume the exogenous restless model and that all arms, when engaged,
are modeled as finite state, irreducible, aperiodic, and reversible
Markov chains. For any increasing sequence $f(t)$ ($f(t) \rightarrow
\infty $ as $t \rightarrow \infty$), if $L(t)$ and $D(t)$ are chosen
such that $L(t) \to \infty$ as $t \to \infty$ ,$\frac{f(t)}{D(t)}
\to \infty $ as $t \to \infty$, and $\frac{D(t)}{L(t)} \to \infty$
as $t \to \infty$, then we have
\begin{eqnarray}
r_{\Phi}(t) \sim o(f(t)\log(t)).
\end{eqnarray}
\end{theorem}
We point out that the conclusion in Theorem~$2$ still holds for the
endogenous restless model, though the proof needs to be modified.
\begin{proof}
See Appendix~B for details.

\end{proof}
\section{Conclusion} \label{sec:conclusion}

In this paper, we studied the decentralized restless multi-armed
bandit problems, where distributed players aim to accrue the maximum
long-term reward without knowing the system reward statistics. Under
the exogenous model where the arm reward status remains static when
not engaged, we proposed a policy to achieve the optimal logarithmic
order of the system regret. Under the endogenous model where the arm
reward status evolves according to an arbitrary random process when
not engaged, we showed that the proposed policy achieves a
logarithmic (weak) regret. Furthermore, we showed that the proposed
policy achieves a complete decentralization where no pre-agreement
or global synchronization among players is required.

\section*{Appendix A. Proof of Theorem~\ref{thm:singleuser}}\label{sec:AppendixA}
We first rewrite the definition of regret as
\begin{eqnarray}\label{eqn:regreform}
r_{\Phi}(t) & = & t\sum_{i=1}^{M}\mu_{\sigma(i)}  - \mbbE_{\Phi}R(t) \nonumber \\
&= &\sum_{i=1}^{N} [\mu_i \mbbE[T_i(t)] -
\mbbE[\sum_{n=1}^{T_i(t)}s_i(t_i(n))]] \nonumber \\
& & +\mbbE[\sum_{n=1}^{T_i(t)}s_i(t_i(n))]] - \mbbE_{\Phi}R(t)
\nonumber \\ & & + t\sum_{i=1}^{M}\mu_{\sigma(i)}- \sum_{i=1}^{N}
\mu_i \mbbE[T_i(t)].
\end{eqnarray}

To bound the first term in (\ref{eqn:regreform}), Lemma~$1$ is
introduced below:

 {\em Lemma $1$}~\cite{Anantharam:87-2}: Let
$Y_1,Y_2, \cdots$ be Markovian with state space $\mathcal{S}$,
matrix of transition probabilities $P$, an initial distribution
$\vec{q}$, and stationary distribution $\vec{\pi}$ ($\pi_s$ is the
stationary probability of state $s$). Let $F_t$ be the
$\sigma$-algebra generated by $Y_1,Y_2, \cdots, Y_t$ and $G$ an
$\sigma$-algebra independent of $Y_{\infty} = \lor Y_t$. Let $T$ be
a stopping time of $\{F_t \lor G\}$. The state (reward) at time $t$
is denoted by $s(t)$. Let $\mu$ denote the mean reward. For any
stopping time $T$, there exists a value $A_P \leq (\min_{s\in
\mathcal{S}} \pi_s)^{-1} \sum_{s\in \mathcal{S}} s $ such that
$\mbbE[\sum_{t=1}^{T}s(t) -\mu T] \leq A_P.$

Using Lemma~$1$ the first term in (\ref{eqn:regreform}) can be
bounded by the following constant:

\begin{eqnarray}\label{eqn:constant}
\sum_{i=1}^{N}[(\min_{s\in \mathcal{S}_i} \pi_s)^{-1}  \sum_{s\in
\mathcal{S}_i} s]
\end{eqnarray}

To show that the regret has a logarithmic order, it is sufficient to
show that the second term plus the third term in
(\ref{eqn:regreform}) has a logarithmic order. These two terms can
be understood as regret caused by two reasons. The first one is
engaging bad arms in the exploration epochs. The second one is not
playing the expected arms in the exploitation epochs. To show the
second term in (\ref{eqn:regreform}) has a logarithmic order, it is
sufficient to show that the regret caused by the two reasons above
have logarithmic orders.

Let $\mbbE[T_{O}(t)]$ denote the time spent on each arm in the
exploration epochs by time $t$ and an upper bound on $T_{O}(t)]$ is:
\begin{eqnarray}\label{eqn:exploration1}
T_{O}(t) \leq \frac{1}{3}[4(3D \ln t +1)-1].
\end{eqnarray}

Consequently the regret caused by engaging bad arms in the
exploration epochs by time $t$ is upper bounded by

\begin{eqnarray}\label{eqn:exploration}
\frac{1}{3}[4(3D \ln t +1)-1]  \left( \sum_{i=1}^{M}\mu_{\sigma(i)}
- \frac{M}{N} \sum_{i=1}^{N}\mu_{\sigma(i)}\right).
\end{eqnarray}

The second reason for regret in the second term of
(\ref{eqn:regreform}) is not playing the expected arms in the
exploitation epochs. Let $t_n$ denote the beginning point to the
$n$th exploitation epoch. Let $\Pr[i,j,n]$ denote the possibility
that arm $i$ has a higher index than arm $j$ at $t_n$, where $\mu_i
< \mu_j$ and $\mu_j \geq \mu_{\sigma(M)}$. It can be shown that:
 \begin{eqnarray}
\Pr[i,j,n] \leq \frac{|\mathcal{S}_i| +
|\mathcal{S}_j|}{\pi_{\min}}(1+ \frac{\epsilon_{\max}
\sqrt{L}}{10s_{\min}}) t_n ^{-1}
\end{eqnarray}

Since different subepochs in the exploitation epochs are symmetric,
the regret in different subepochs are the same. In the first
subepoch, player $k$ aims at arm $\sigma(k)$. In the model where
players in conflict share the reward, player $k$ failing to identify
arm $\sigma(k)$ in the first subepoch of the $n$th exploitation
epoch can lead to a regret no more than $\mu_{\sigma(k)} 2 \times
4^{n-1} $. In calculating the upper bound for regret, for player
$M$, we can assume that playing the arm $\sigma(M+1)$ to arm
$\sigma(N)$ can contribute to the total reward. Thus an upper bound
for regret in the $n$th exploitation epoch can be obtained as
\begin{eqnarray}
&&2M 4^{n-1}(1+ \frac{\epsilon_{\max} \sqrt{L}}{10s_{\min}}) t_n
^{-1}[ \sum_{i=1}^{M-1} \sum_{j=1,j \neq i}^N
\mu_i\frac{|\mathcal{S}_i| + |\mathcal{S}_j|}{\pi_{\min}}\nonumber  \\
&&  + \sum_{j=1}^{M-1} \mu_M \frac{|\mathcal{S}_M| +
|\mathcal{S}_j|}{\pi_{\min}} \nonumber  \\ &&+ \sum_{j=M+1}^N
(\mu_M-\mu_j)\frac{|\mathcal{S}_i| + |\mathcal{S}_j|}{\pi_{\min}} ]
\end{eqnarray}

By time $t$, at most $(t-N)$ time slots have been spent on the
exploitation epochs. Thus
\begin{eqnarray}\label{eqn:ubexploitation}
n_I(t) \leq \lceil \log_4 (\frac{3}{2} (t- N) +1) \rceil .
\end{eqnarray}

From the upper bound on the number of the exploitation epochs given
in (\ref{eqn:ubexploitation}), and also the fact that $t_n \geq
\frac{2}{3} 4^{n-1}$, we have the following upper bound on regret
caused in the exploitation epochs by time $t$ (Denoted by
$r_{\Phi,I}(t)$):

\begin{eqnarray} \label{eqn:UPI}
r_{\Phi,I}(t)& \leq & 3\lceil \log_4 (\frac{3}{2} (t- N) +1)
\rceil(1+ \frac{\epsilon_{\max} \sqrt{L}}{10s_{\min}}) \nonumber \\
 &  & \qquad \sum_{i=1}^{M-1} \sum_{j=1,j \neq i}^N
\mu_{\sigma(i)}\frac{|\mathcal{S}_{\sigma(i)}| + |\mathcal{S}_{\sigma(j)}|}{\pi_{\min}} \nonumber \\
&  & +  3 \lceil \log_4 (\frac{3}{2} (t- N) +1) \rceil(1+
\frac{\epsilon_{\max} \sqrt{L}}{10s_{\min}}) \nonumber \\
 &  & \qquad \sum_{j=M+1}^N
(\mu_{\sigma(M)}-\mu_{\sigma(j)})\frac{|\mathcal{S}_{\sigma(i)}| +
|\mathcal{S}_{\sigma(j)}|}{\pi_{\min}} \nonumber \\&  & + 3 \lceil
\log_4 (\frac{3}{2} (t- N) +1) \rceil(1+ \frac{\epsilon_{\max}
\sqrt{L}}{10s_{\min}}) \nonumber \\
 &  & \qquad  \sum_{j=1}^{M-1} \mu_{\sigma(M)}
\frac{|\mathcal{S}_{\sigma(M)}| +
|\mathcal{S}_{\sigma(j)}|}{\pi_{\min}}
\end{eqnarray}

Combining~(\ref{eqn:regreform})~(\ref{eqn:constant})~(\ref{eqn:exploration})~(\ref{eqn:UPI}),
we can get the upper bound of regret:
\begin{eqnarray}
r_{\Phi}(t)& \leq &\frac{1}{3}[4(3D \ln t +1)-1]  \left(
\sum_{i=1}^{M}\mu_{\sigma(i)}
- \frac{M}{N} \sum_{i=1}^{N}\mu_{\sigma(i)}\right)\nonumber \\
 &  & + 3\lceil \log_4 (\frac{3}{2} (t- N) +1)
\rceil(1+ \frac{\epsilon_{\max} \sqrt{L}}{10s_{\min}}) \nonumber \\
 &  & \qquad \sum_{i=1}^{M-1} \sum_{j=1,j \neq i}^N
\mu_{\sigma(i)}\frac{|\mathcal{S}_{\sigma(i)}| + |\mathcal{S}_{\sigma(j)}|}{\pi_{\min}} \nonumber \\
&  & +  3 \lceil \log_4 (\frac{3}{2} (t- N) +1) \rceil(1+
\frac{\epsilon_{\max} \sqrt{L}}{10s_{\min}}) \nonumber \\
 &  & \qquad \sum_{j=M+1}^N
(\mu_{\sigma(M)}-\mu_{\sigma(j)})\frac{|\mathcal{S}_{\sigma(i)}| +
|\mathcal{S}_{\sigma(j)}|}{\pi_{\min}} \nonumber \\&  & + 3 \lceil
\log_4 (\frac{3}{2} (t- N) +1) \rceil(1+ \frac{\epsilon_{\max}
\sqrt{L}}{10s_{\min}}) \nonumber \\
 &  & \qquad  \sum_{j=1}^{M-1} \mu_{\sigma(M)}
\frac{|\mathcal{S}_{\sigma(M)}| +
|\mathcal{S}_{\sigma(j)}|}{\pi_{\min}}
\nonumber \\
&&+ \sum_{i=1}^{N}[(\min_{s\in \mathcal{S}_i} \pi_s)^{-1} \sum_{s\in
\mathcal{S}_i} s]
\end{eqnarray}

Next we consider the model where no player in conflict gets reward.
In the first subepoch of the $n$th exploitation epcoh, each mistake
by player $k$ can cause regret more than $\mu_{\sigma(k)} 2 \times
4^{n-1}$. Assuming each mistake can cause
$\sum_{i=1}^{M}\mu_{\sigma(i)} 2 \times 4^{n-1}$ regret leads to the
following upper bound for regret under this conflict model:

\begin{eqnarray}
r_{\Phi}(t)& \leq &3 \lceil \log_4 (\frac{3}{2} (t- N) +1) \rceil(1+
\frac{\epsilon_{\max}
\sqrt{L}}{10s_{\min}})\nonumber \\
 &  & \qquad  (\sum_{i=1}^{M}\mu_{\sigma(i)}) \sum_{i=1}^{M}
\sum_{j=1,j \neq i}^N \frac{|\mathcal{S}_{\sigma(i)}| +
|\mathcal{S}_{\sigma(i)}|}{\pi_{\min}}\nonumber \\
 &  & +
\frac{1}{3}[4(3D \ln t +1)-1] \nonumber \\
 &  & \qquad \left( \sum_{i=1}^{M}\mu_{\sigma(i)}
- \frac{M}{N} \sum_{i=1}^{N}\mu_{\sigma(i)}\right)
\nonumber \\
&&+ \sum_{i=1}^{N}[(\min_{s\in \mathcal{S}_i} \pi_s)^{-1} \sum_{s\in
\mathcal{S}_i} s]
\end{eqnarray}

\section*{Appendix B. Proof of Theorem~\ref{thm:arbitraryclose}}\label{sec:AppendixB}

The choice of $L(t)$ and $D(t)$ implies that $D(t) \to \infty$ as $t
\to \infty$. The regret has three parts: the transient effect of
arms, the regret caused by playing bad arms in the exploration
epochs, and the regret caused by mistakes in the exploitation
epochs. It will be shown that each part part of the regret is on a
lower order than $f(t)\log(t)$. The transient effect of arms is the
same as in Theorem~$1$. Thus it is upper bounded by a constant
independent of time $t$ and is on a lower order than $f(t)\log(t)$.

The regret caused by playing bad arms in the exploration epochs is
bounded by
\begin{eqnarray}\label{eqn:1}
\frac{1}{3}[4(3D(t) \ln t +1)-1]  \left(
\sum_{i=1}^{M}\mu_{\sigma(i)} - \frac{M}{N}
\sum_{i=1}^{N}\mu_{\sigma(i)}\right).
\end{eqnarray}
Since $\frac{f(t)}{D(t)} \to \infty $ as $t \to \infty$, the part of
regret in~(\ref{eqn:1}) is on a lower order than $f(t)\log(t)$.

For the regret caused by playing bad arms in the exploitation
epochs, it is shown below that the time spent on a bad arm $i$ can
be bounded by a constant independent of $t$.

Since $\frac{D(t)}{L(t)} \to \infty $ as $t \to \infty$, there
exists a time $t_1$ such that $\forall t \geq t_1$, $D(t) \geq
\frac{4L(t)}{(\min_{j\leq M}( \mu_{\sigma(j)} -
\mu_{\sigma(j+1)}))^2}$. There also exists a time $t_2$ such that
$\forall t \geq t_2$, $L(t) \geq \frac{1}{\epsilon_{\min}}(7 \frac
{20s^2_{\max} |\mathcal{S}|_{\max}^2}{(3-2\sqrt{2})}+ 10
s^2_{\max})$. The time spent on playing bad arms before $t_3 =
\max(t_1,t_2)$ is at most $t_3$, and the caused regret is at most
$(\sum_{j=1}^M\mu_{\sigma(j)})t_3$. The regret caused by mistakes
after $t_3$ is upper bounded by $6(1+ \frac{\epsilon_{\max}
\sqrt{L}}{10s_{\min}})(\sum_{i=1}^{M}\mu_{\sigma(i)}) \sum_{i=1}^{M}
\sum_{j=1,j \neq i}^N \frac{|\mathcal{S}_{\sigma(i)}| +
|\mathcal{S}_{\sigma(i)}|}{\pi_{\min}}$. Thus the regret caused by
mistakes in the exploitation epochs is on a lower order than
$f(t)\log(t)$.

Because each part of the regret is on a lower order than
$f(t)\log(t)$, the total regret is also on a lower order than
$f(t)\log(t)$.

\end{document}

%% file: macros.tex
\def\boxit#1{\vbox{\hrule\hbox{\vrule\kern3pt
        \vbox{\kern3pt#1\kern3pt}\kern3pt\vrule}\hrule}}

\def\reals{ { {\rm  I \kern-0.15em R }  } }
\def\complex{ {\,{{\rm C} \kern-0.50em \raise0.20ex {  |}}\, }}

\def\Rbf{{\bf R}}

\def\be{\begin{equation}}
\def\ee{\end{equation}}

\def\defeq{{\stackrel{\Delta}{=}}}
\def\scalefig#1{\epsfxsize #1\textwidth}
\def\Rxx{\Rbf_{\ssstyle X\kern-.1em X}}

\let\ssstyle=\scriptscriptstyle

\def\etal{{\it et al. \/}}

\def\ie{{\it i.e.,\ \/}}
\def\Kout{\setbox1=\hbox{\Huge\bf K}\hbox to
1.05\wd1{\hspace{.05\wd1}
}}
\def\Sout{\setbox1=\hbox{\Huge\bf S}\hbox to 1.05\wd1{\hspace{.05\wd1}
}}